\documentclass[11pt,leqno]{article} 
\usepackage{amsmath,amsfonts,amscd,ProtocoloEulerIH,url}

\title{\Large\bf  {Intersection cohomology of the circle actions}\footnote{This work 
has been partially supported by the UPV-EHU
127.310-E-14790/2002 (second author) and the 
EC0S-Nord Projet V00M01.}}

\author{
Gabriel Padilla\thanks{Escuela de Matemática. Universidad Central
Venezuela. Caracas 1010-Venezuela. {\sl  gpadilla@euler.ciens.ucv.ve}. }
\\ {\small Universidad Central de Venezuela}
\and  
Martintxo Saralegi-Aranguren\thanks{Fédération CNRS Nord-Pas-de-Calais FR 2956.
UPRES-EA 2462 LML. 
Faculté Jean Perrin.
Université d'Artois.   Rue Jean Souvraz SP 18.   62 307 Lens Cedex - 
France.   
{\sl saralegi@euler.univ-artois.fr}. }
\\ {\small Université d'Artois }
}

\begin{document}  
\maketitle

\begin{abstract}
A classical result says that a free action of the circle $\S^1$ on a
topological space $X$  is geometrically classified by the orbit
space $B$ and by a cohomological class $e \in \coho{H}{2}{B,\Z}$, the
Euler class.  When the action is not free  we have a difficult open
question:
\begin{center}
$\Pi$ : ``Is the space $X$ determined by the orbit space $B$ and the
Euler class?''
\end{center}

The main result of this work is a step towards the understanding of the above
question in the category of unfolded pseudomanifolds.  We prove
that the orbit space $B$ and the Euler class determine:
\begin{itemize}
\item the intersection cohomology of $X$, 
\item the real homotopy
type of $X$.
\end{itemize}
 \end{abstract}

In this work, we give an answer to the question $\Pi$ in the category of unfolded 
pseudomanifolds. The object studied  are the modelled actions
$\Phi \colon \sbat \times X \to X$. Here, the total space $X$ is an 
unfolded pseudomanifold and the action $\Phi$ preserves this 
structure in such a way that the orbit space $B$ is still an unfolded 
pseudomanifold.

A priori, the action $\Phi$ classifies the strata of $X$ in two 
types: the mobile strata (containing one-dimensional orbits), and the
fixed strata (containing the fixed points).  But we see in this work
that we need a finer classification: a fixed stratum $S$ can be
perverse or not perverse.  The stratum $S$ is perverse when the action of $\S^1$ on its link
is \color{blue} not \normalcolor cohomologically trivial.

\color{blue}In the context of singular actions, the meaning of ``Euler
class'' it is not clear: there are non trivial circle actions having a contractible orbit space $B$.
This Euler class $e$ can be recovered by using 
the de Rham intersection cohomology $\coho{\IH}{*}{-}$.  It has been proved that $e$ lives in
$\lau{\IH}{2}{\per{e}}{B}$ where the Euler perversity $\per{e}$ takes
the following values
$$
    \per{e}(S) = 
    \left\{
    \begin{array}{ll}
    0 & \hbox{when $S$ mobile stratum}, \\
1 & \hbox{when $S$ not perverse fixed stratum}\\
2 & \hbox{when $S$ perverse stratum}
    \end{array}
    \right.
    $$
 (cf. \cite[5.7]{Pa}). Notice
    that the Euler class contains the geometrical information about the
    nature of the strata.\normalcolor

\medskip

The main result of this work is the following: the orbit space $B$
of a modelled action and the Euler
class $e \in \lau{\IH}{2}{\per{e}}{B}$ determine the intersection
cohomology of $X$ (cf.  Corollary \ref{uno}), the real homotopy type
of $X$ (cf.  Corollary \ref{dos}) and the perverse real homotopy type of $X$
(cf.  Corollary \ref{tres}).  The main tool we use is the Gysin
sequence constructed for $\Phi$ in \cite{Pa}.

\vspace{,5cm}

%

\section{Intersection cohomology of modelled actions}

  We recall in this Section the main results of \cite{Pa} we are going to use in this work.  
   
  \prg {\bf Modelled actions.} A reasonable action of the circle on a 
  stratified pseudomanifold must produce a stratified pseudomanifold as 
  orbit space. These are the $\sbat$-pseudomanifolds of \cite[Section 
  4]{Po}. In this 
  work we shall use a variant of this concept, the modelled action $\Phi \colon \sbat \times X \to X$ of the circle $\sbat$ on an unfolded pseudomanifold $X$, since the unfolded pseudomanifolds support the  (de Rham) intersection cohomology (cf. \cite{Pa}).
   We list  below the main properties   of a modelled action
  $\Phi \colon \sbat \times X \TO X$ of the circle  $ \sbat$ on an unfolded 
  pseudomanifold $X$. We denote by $B =X/\sbat$ the orbit space and by $\pi \colon X \to B$ the canonical projection.

  \Am {\em The isotropy subgroup $\sbat_x$ is the same for each $x \in S$.  
  It will be denoted by $\sbat_S$}.

  \am {\em For each regular stratum $R$ we have $\sbat_{R} = \{ 1 \}$}.

  \am {\em For each  singular stratum $S$  with $\sbat_{S}= \sbat$, the 
  action $\Phi$ induces a modelled action $\Phi_{L_S} 
  \colon 
  \sbat \times L_S \to L_S$, where $L_S$ is the link of $S$}.

  \am {\em The orbit space $B$ is an unfolded 
  pseudomanifold, 
  relatively to the stratification ${\cal S}_B = \{ \pi(S) \ / \ S \in {\cal 
  S}_X\}$, and the projection  $\pi \colon X \to B $ is an unfolded morphism}.

  \am {\em The assignment $S \mapsto \pi (S)$ induces  the bijection $\pi_{\mathcal S}   \colon {\cal S}_X \to 
  {\cal 
  S}_B$.}

  \bigskip

  The action $\Phi$ classifies 
    the  strata of $X$
  in two types: 
  the stratum $S$ is  {\em mobile} when $\sbat_S$ is finite
  and it is  {\em fixed} when $\sbat_S = 
    \sbat$.   In this work, we need another classification for the 
    fixed strata. A fixed stratum $S$ is {\em perverse} when 
    $\coho{H}{*}{L_S \menos \Sigma_{L_S}}\color{blue}\not= \normalcolor
    \coho{H}{*}{(L_S \menos \Sigma_{L_S})/\sbat} 
  \otimes 
  \coho{H}{*}{\sbat}$, where $\Sigma_{L_S}$ is the singular part of the link $L_{S}$ \color{blue}(cf. \cite[5.6 (3)]{Pa})\normalcolor.

  \prg {\bf Examples.}
  Consider $B = c\S^2$. Essentially, there are three different modelled 
  actions having $B$ as the orbit space.

  $
  \begin{array}{lcl}
  \Phi_1 \colon \sbat \times c\S^3 \TO c\S^3
  &\hbox{ defined by } & 
  \Phi_1(z,[(u,v),t]) = [(z\cdot u, z \cdot v),t],
  \\[,3cm]
  \Phi_2 \colon \sbat \times c(\S^2 \times\S^1 )\TO c(\S^2 \times\S^1 )
  &\hbox{ defined by } &  
  \Phi_2(z,[(x,w),t]) = [(x, z \cdot w),t], \hbox{ and }
  \\[,3cm]
  \Phi_3 \colon \sbat \times  \left( c(\S^2) \times \S^1 \right)\TO 
  \left( c(\S^2) \times\S^1 \right)
  &\hbox{ defined by } &  
  \Phi_3(z,([x,t],w) = ([x, t] ,z \cdot w).
  \end{array}
  $

  \medskip

The 
  difference between these actions lies on the geometrical nature 
  of the singular stratum $\{ \vartheta \}$ (vertex)  of $B$. In fact, in the first case the stratum 
  $\{ \vartheta \}$ comes from  a perverse stratum,
  in the second case   the stratum   $\{ \vartheta \}$  comes from   a 
   non-perverse  fixed stratum and in the third case 
   the stratum 
   $\{ \vartheta \}$ comes from   a mobile stratum.

\prg {\bf Gysin sequence}.  Since the Lie group $\sbat$ is connected and compact, the subcomplex 
 of the invariant perverse forms 
  computes the 
 intersection cohomology\footnote{For the notions related with the intersection cohomology, we refer the reader to \cite[Section 3]{S2}.} of $X$. In fact,  for any perversity $\per{p}$, the inclusion
 $\left( \lau{\Omega}{*}{\per{p}} {X}\right)^{\sbat} \hookrightarrow
 \lau{\Omega}{*}{\per{p}} {X}$ induces an isomorphism in cohomology.
 This complex can described in terms of basic data as follows.
 Consider
 the graded  
 complex
 \be
\label{inv}
\lau{I\Omega}{*}{\per{p}} {X}= \left\{ (\alpha,\beta) \in
\lau{\Pi}{*}{}{B} \oplus \lau{\Om}{*-1}{\per{p}-\per{x}}{B}  \Big/
 \left[
 \begin{array}{l}
 ||\alpha||_{\pi(S)} \leq \per{p}(S) \\[,3cm] ||d\alpha + (-1)^{|\beta|}
 \beta \wedge \epsilon||_{\pi(S)} \leq \per{p}(S)
 \end{array}
 \right]
 \hbox{ if } S \in \mathcal{S}^{^{sing}}_{X} \right\} \ee endowed with
 the differential $ D(\alpha,\beta) = (d\alpha + (-1)^{|\beta|} \beta
 \wedge \epsilon , d\beta) $.  Here $| -|$ stands for the degree of the
 form, $\epsilon \in \coho{\Pi}{2}{B}$ is an Euler form and $\per{x}$ is the characteristic perversity defined by 
 $
\per{x}(\pi(S)) = \left\{
\begin{array}{ll}
1& \hbox{if $S$ is a fixed stratum }\\
0 & \hbox{if $S$ is a mobile stratum} 
\end{array}
\right.
$. The assignment $(\alpha,\beta) \mapsto \pi^{*}\alpha +
 \pi^{*}\beta \wedge \chi$ establishes a differential graded
 isomorphism between $\lau{I\Om}{}{\per{p}}{X}$ and $ \left(
 \lau{\Omega}{*}{\per{p}} {X}\right)^{\sbat}$.
 
From \refp{inv} we have the short
exact sequence
$$
0
\TO
\lau{\Om}{*}{\per{p}}{B} 
\stackrel{\ib{\pi}{\per{p}}}{\TO}
\lau{I\Om}{*}{\per{p}}{X} 
\stackrel{\ib{\oint}{\per{p}}}{\TO}
\lau{{\cal G}}{\color{blue}*-1}{\per{p}}{B}
\TO
0,
$$
where
\begin{itemize} 

\item The {\em Gysin term} $\lau{{\cal G}}{\color{blue}*-1}{\per{p}}{B}$ is the 
differential 
complex
 $$ 
\left\{
\beta \in \lau{\Om}{*-1}{\per{p} - \per{x}}{B} \ \Big/ \ \exists \alpha \in 
\lau{\Pi}{*}{}{B} \hbox{ with } 
\left[
\begin{array}{l}
||\alpha||_{\pi(S)} \leq \per{p}(S)\hbox { and } \\[,3cm] ||d\alpha +
(-1)^{|\beta|}\beta \wedge \epsilon||_{\pi(S)} \leq \per{p}(S)
\end{array}
\right]
\hbox{ if } S \in {\cal S}^{^{sing}}_X
\right\},
$$

\item $\ib{\oint}{\per{p}} (\alpha, \beta) = \beta$, and

\item $\ib{\pi}{\per{p}}(\om) = \pi^{*}\om$.

\end{itemize}
The associated long exact sequence
 $$
 \cdots 
 \TO
 \lau{\IH}{i+1}{\per{p}}{X} \stackrel{\ib{\oint}{\per{p}}}{\TO}
 \coho{H}{i}{\lau{{\cal G}}{*}{\per{p}}{B}}  
 \stackrel{\ib{\eub}{\per{p}}}{\TO}
 \lau{\IH}{i+2}{\per{p}}{B} 
 \stackrel{\ib{\pi}{\per{p}}}{\TO}
 \lau{\IH}{i+2}{\per{p}}{X} 
 \TO
 \dots,
 $$
where
$\ib{\eub}{\per{p}} [\beta] = [d\alpha + (-1)^{|\beta |}\beta \wedge 
 \epsilon]$, is the {\em Gysin sequence}.

  Recall that the Euler perversity $\per{e}$ is defined by
  $
    \per{e}(S) = 
    \left\{
    \begin{array}{ll}
    0 & \hbox{when $S$ mobile stratum}, \\
1 & \hbox{when $S$ not perverse fixed stratum}\\
2 & \hbox{when $S$ perverse stratum}
    \end{array}
    \right.
    $
    So, the Euler class $e =  [\epsilon]$ belongs to $\lau{\IH}{2}{\per{e}}{B}$. This class detects the perverse strata: a fixed stratum is perverse iff the Euler class $e_S \in \lau{\IH}{2}{\per{e}}{L_S/\sbat}$ of the action $\Phi_{L_S} \colon \sbat \times L_S \to M_S$, \color{blue} does not vanish \normalcolor (see (MA.iii)).
In the next Section, we shall use the following Lemma
  \bL
  Let $\per{p}$ be a perversity with $\per{p} \geq \per{e}$. If $X$ is connected and normal, then
    \be
 \label{g}
 \coho{H}{0}{\lau{{\cal 
 G}}{*}{\per{p}}{B}} \cong \R \ \ \hbox{ and }\ib{\eub}{\per{p}}(1) = e.
 \ee
 \eL
 \pro  Condition
 $\per{p} \geq \per{e}$ implies $1 \in\lau{{\cal G}}{*}{\per{p}}{B}$. 
 Since $X$ is connected and normal, then  the regular part $B \menos\Sigma_{B}$ is connected.  Then
 $ \coho{H}{0}{\lau{{\cal G}}{*}{\per{p}}{B}} \cong \R$.  Finally, the
 definition of $\ib{\eub}{\per{p}}$ gives $\ib{\eub}{\per{p}}(1) = [\epsilon ] =
 e$.
\qed

\section{Perverse algebras}

Although the intersection cohomology $\lau{\IH}{*}{\per{p}}{X}$ is not an algebra, we recover this structure  by considering all the perversities together. These are the perverse algebras we present in this Section.

\prg {\bf Perverse algebras}.  
A {\em perverse set} is a triple $(\mathcal{P},+,\leq)$ where 
$(\mathcal{P},+)$ is an abelian semi-group with an unity element 
$\per{0}$ and $(\mathcal{P},\leq)$ is a poset verifying the 
compatibility condition:
$$
\per{p}\leq \per{q} \hbox{ and } \per{p}'\leq \per{q}' \Longrightarrow 
\per{p} + \per{p}' \leq \per{q}' + \per{q}',
\ \ \hbox{ for } \per{p}, \per{q}, \per{p}', \per{q}' \in \mathcal{P}.
$$
In order to simplify the writing, we shall say that $\mathcal{P}$ is a 
perverse set.

 A {\em dgc perverse algebra} (or simply a perverse algebra) is a quadruple
 $\mbox{\boldmath $E$} = (E,\iota,\wedge,d)$ where
 \begin{itemize}
     \item[-] $E =\displaystyle{\bigoplus_{\per{p} \in \mathcal{P}}}\ib{E}{\per{p}}$
     where each $\ib{E}{\per{p}}$ is a graded (over $\Z$) vector space,

    \item[-]
    $\iota = \left\{ \ib{\iota}{\per{p},\per{q}} \colon \ib{E}{\per{p}} \to
    \ib{E}{\per{q}} \ / \ \per{p} \leq \per{q} \right\}$ is a family of
 graded linear morphisms,    and
    
  \item[-] $(E,d,\wedge)$ is a dgc algebra, 
    \end{itemize}
    verifying
$$
\begin{array}{lll}
+ \ \ib{\iota}{\per{p},\per{p}} = \Ide \hspace{1cm}&

+ \  \ib{\iota}{\per{q} ,\per{r}} \rondp \ib{\iota}{\per{p} ,\per{q}} =
      \ib{\iota}{\per{p}, \per{r}} & 
      
+ \ \wedge \left( \ib{E}{\per{p}} \times
	   \ib{E}{\per{p}'}\right) \subset \ib{E}{\per{p} +\per{p}'}\\[,3cm]
	   
      + \ d\left( \ib{E}{\per{p}}\right) \subset \ib{E}{\per{p}}
&
 
  + \ \ib{\iota}{\per{p}+\per{p}',\per{q}+\per{q}'}(a \wedge
      a') = \ib{\iota}{\per{p},\per{q}}(a) \wedge
      \ib{\iota}{\per{p}',\per{q}'}(a')  \hspace{1cm}
&
+ \ d \rondp \ib{\iota}{\per{p},\per{q}} = 
      \ib{\iota}{\per{p},\per{q}} \rondp d  \\
\end{array}
$$
 \nt Here, $\per{p} \leq \per{q} \leq \per{r}$, $\per{p}' \leq \per{q}'$, $a \in
 \ib{E}{\per{p}}$ and $a'\in \ib{E}{\per{p}'}.
$ \smallskip

 Associated to a dgc perverse algebra $\mbox{\boldmath $E$} =
 (E,\iota,\wedge,d)$ we have another dgc perverse algebra, namely, its
 cohomology $\lau{\mbox{\boldmath $H$} }{}{}{\mbox{\boldmath $E$}}=
 \left(\displaystyle{\bigoplus_{\per{p} \in
 \mathcal{P}}}\lau{H}{}{}{\ib{E}{\per{p}},d},\iota,\wedge,0\right)$,
 where $\iota$ and $\wedge$ are induced by the previous $\iota$ and
 $\wedge$.

A {\em dgc perverse morphism} (or simply {\em perverse morphism})
{\boldmath $f$} between two perverse algebras $\mbox{\boldmath $E$} =
(E,\iota,\wedge,d)$ and $\mbox{\boldmath $E'$} =
(E',\iota',\wedge',d')$ is given by a family $ \mbox{\boldmath $f$} =
\left\{ \ib{f}{\per{p}} \colon \ib{E}{\per{p}} \to \ib{E}{\per{p}}
\right\}$ of differential graded morphisms verifying
\begin{equation}
\label{bat}
 \ib{\iota'}{\per{p},\per{q}} \rondp \ib{f}{\per{p}}=
	\ib{f}{\per{q}} \rondp \ib{\iota}{\per{p},\per{q}}
	\end{equation}
	and
	\begin{equation}
	\label{bi}
	\ib{f}{\per{p}+\per{p}'} (a \wedge b) = \ib{f}{\per{p}}(a) \wedge
	\ib{f}{\per{p}'}(b).
\end{equation}
 Here, $\per{p}\leq \per{q}$, $a \in
 \ib{E}{\per{p}}$ and $b \in \ib{E}{\per{p}'}$.  We shall denote the
 perverse morphism by $\mbox{\boldmath $f$} \colon \mbox{\boldmath $E$} \to
 \mbox{\boldmath $E'$}$.  It induces the  perverse morphism 
 $\mbox{\boldmath $f$}\colon \lau{\mbox{\boldmath $H$} }{}{}{\mbox{\boldmath $E$}}
 \to
 \lau{\mbox{\boldmath $H$} }{}{}{\mbox{\boldmath $E$}'}$, 
 defined by $\ib{f}{\per{p}}[a] = [\ib{f}{\per{p}}(a)]$ 
 for each $\per{p} $ and $[a] \in 
 \lau{H}{}{}{\ib{E}{\per{p}},d}$.
 
 When each $\ib{f}{\per{p}}$ is an isomorphism, we
 shall say that $\mbox{\boldmath $f$}$ is a {\em dgc perverse
 isomorphism} (or simply {\em perverse isomorphism}).
 It induces the perverse isomorphism $\mbox{\boldmath $f$} \colon
 \lau{\mbox{\boldmath $H$} }{}{}{\mbox{\boldmath $E$}} \to
 \lau{\mbox{\boldmath $H$} }{}{}{\mbox{\boldmath $E$}'}$.

\prg {\bf Perverse algebras and modelled actions}. Fix $\Phi \colon \sbat \times X \to X$ a modelled action. The family of perversities $\mathcal{P}_X$ of $X$ has a partial
order $\leq$ and an abelian law $+$ in such a way that $\mathcal{P}_{X}$
is a perverse set.  In the same way, $\mathcal{P}_{B}$
is a perverse set.
Since the two posets ${\cal S}^{^{sing}}_B$
and ${\cal S}^{^{sing}}_X $ are isomorphic (cf. (MA.v)), then the perverse sets $\mathcal{P}_{B}$ and
$\mathcal{P}_{X}$ are isomorphic through the map $\per{p} \mapsto
\per{p}\rondp \pi$ (cf. (MA.iv)).  In the sequel, we shall
identify these two perverse sets.

Associated to the modelled action $\Phi$, we have the following dgc perverse algebras.
\begin{itemize}
\item[+]  The {\em perverse de Rham algebra}: 
$
\lau{\mbox{\boldmath $\Om$}}{}{}{X} = \left( {\displaystyle
\lau{\Om}{}{}{X} = \bigoplus_{\per{p} \in \mathcal{P}_{X}}}
\lau{\Om}{}{\per{p}}{X},\iota,\wedge,d \right).
$
\item[+] The {\em intersection cohomology algebra}:
$
\lau{\mbox{\boldmath $\IH$}}{}{}{X} = \left( {\displaystyle
\lau{\IH}{}{}{X} = \bigoplus_{\per{p} \in \mathcal{P}_{X}}}
\lau{\IH}{}{\per{p}}{X},\iota,\wedge,0\right).
$
\end{itemize}
Analogously for $B$. 
The quadruple
 $
 \lau{\mbox{\boldmath $I\Om$}}{}{}{X} = \left( {\displaystyle
 \lau{I\Om}{}{}{X} = \bigoplus_{\per{p} \in \mathcal{P}_{X}}}
 \lau{I\Om}{}{\per{p}}{X},\iota,\wedge,D\right)
 $
 is a also perverse algebra.  Here, the wedge product is defined by $(\alpha,\beta)
 \wedge (\alpha',\beta') = (\alpha \wedge \alpha' ,
 (-1)^{|\alpha'|}\beta\wedge \alpha' + \alpha \wedge \beta')$.
 A
 straightforward calculation shows that the operator
 \begin{equation}
     \label{delta}
 {\mbox{\boldmath $\Delta$} } = \{ \ib{\Delta}{\per{p}}\}\colon
 \lau{\mbox{\boldmath $I\Om$}}{}{}{X} \to \lau{\mbox{\boldmath $\Om$}}{}{}{X},
 \end{equation}
 defined by
 $
 \ib{\Delta}{\per{p}}(\alpha, \beta) = \pi^{*}\alpha + 
 \pi^{*}\beta \wedge \chi, $ induces a perverse isomorphism in cohomology.

 For each perversity $\per{p}$ we have the linear morphism
 $\ib{\rho}{\per{p}} \colon \lau{\Om}{*}{\per{p}}{B} \to \lau{I\Om}{*}{\per{p}}{X}$
 defined by $\ib{\rho}{\per{p}}(\alpha) = (\alpha,0)$.  The operator $
 \mbox{\boldmath $\rho$} = \{ \ib{\rho}{\per{p}}\}\colon
 \lau{\mbox{\boldmath $\Om$}}{}{}{B}\to \lau{\mbox{\boldmath
 $I\Om$}}{}{}{X} $ is a perverse morphism.  It induces the perverse
 morphism $ \mbox{\boldmath $\pi$} = \mbox{\boldmath $\Delta$} \rondp \mbox{\boldmath $\rho$} \colon \lau{\mbox{\boldmath
 $\IH$}}{}{}{B}\to \lau{\mbox{\boldmath $\IH$}}{}{}{X}.$

\section{Cohomological classification of modelled actions}

We considered in this Section a modelled action $\Phi \colon \sbat
\times X \to X$ whose orbit space is a fixed unfolded pseudomanifold $B$. 
We prove that the intersection cohomology
algebra and the (perverse) real homotopy type of $X$ are determined by the
Euler class.

\prg {\bf Fixing the orbit space}. 
Consider $\Phi_{1} \colon \sbat \times X_{1}\to X_{1}$ and $\Phi_{2}
\colon \sbat \times X_{2}\to X_{2}$ two modelled actions and write
$B_{1}$ and $B_{2}$ the two orbit spaces.
Consider $f \colon B_{1} \to B_{2}$ an unfolded isomorphism.  The two posets
${\cal S}^{^{sing}}_{B_{1}}$ and ${\cal S}^{^{sing}}_{B_{2}}$ are
isomorphic through the map $\pi_{1}(S) \mapsto f(\pi_{1}(S))$.  The
perverse sets $\mathcal{P}_{B_{1}}$ and $\mathcal{P}_{B_{2}}$ are
isomorphic through the map $\per{p} \mapsto \per{p}\rondp f^{-1}$.  In the sequel, we shall identify this two perverse sets in
order to compare the perverse de Rham algebras of $X_{1}$ and $X_{2}$.

The induced map $f^{*} \colon \lau{\Pi}{*}{}{B_{2}} \to \lau{\Pi}{*}{}{B_{1}}$
is a well defined differential graded isomorphism.  It preserves the
perverse degree.  For each perversity $\per{p}$ we write
$\ib{f}{\per{p}} \colon \lau{\Om}{*}{\per{p}}{B_{2}} \to
\lau{\Om}{*}{\per{p}}{B_{1}}$ the differential graded isomorphism
defined by $\ib{f}{\per{p}}(\alpha) = f^{*}\alpha$.  The operator $
\mbox{\boldmath $f$} = \{ \ib{f}{\per{p}}\}\colon \lau{\mbox{\boldmath
$\Om$}}{}{}{B_{2}}\to \lau{\mbox{\boldmath $I\Om$}}{}{}{B_{1}}, $ is a
perverse isomorphism.  It induces the perverse isomorphism $
\mbox{\boldmath $f$}\colon \lau{\mbox{\boldmath $\IH$}}{}{}{B_{2}}\to
\lau{\mbox{\boldmath $\IH$}}{}{}{B_{1}}.  $

The unfolded isomorphism $f$ is {\em optimal} when it preserves the nature
of the strata, that is, when it sends the fixed (resp.  perverse,
resp.  non-perverse) strata into fixed (resp.  perverse, resp. 
non-perverse) strata.  In this case, the two Euler perversities are
equal: $\per{e}_{1}(\pi_{1}(S)) = \per{e}_{2}(f(\pi_{1}(S)))$ for each
singular stratum $S \in {\cal S}^{^{sing}}_{X_{1}}$.  We shall write
$\per{e}$ for this Euler perversity.

Now we can compare the two Euler classes
$e_{1}\in \lau{\IH}{2}{\per{e}}{B_{1}}$ and $e_{2} \in \lau{\IH}{2}{\per{e}}{B_{2}}$
We shall say that $e_{1}$ and $e_{2}$ are {\em proportional} if there
exists a number $\lambda \in \R\backslash\{0\}$ such that 
$\ib{f}{\per{e}}(e_{2}) = \lambda \cdot e_{1}$.  As we are going to see, this
is the key test for the comparison between the de Rham algebras of $X_{1}$ ands
$X_{2}$.

\smallskip 

Finally, we say that the actions $\Phi_{1}$ and $\Phi_{2}$ {\em
have a common  orbit space} if there exists an optimal isomorphism between
theirs orbit spaces.

\bigskip

The three main results of this work come from this Proposition.
\bp
\label{comp}
Let $X_{1}$, $X_{2}$ be two connected normal unfolded pseudomanifolds. 
Consider two modelled actions $\Phi_1 \colon \sbat \times X_1 \to X_1$
and $\Phi_2 \colon \sbat \times X _2\to X_2$.  Let us suppose that
there exists an unfolded isomorphism $f \colon B_{1} \to B_{2}$
between the associated orbit spaces.  Then, the two following
statements are equivalent:

\Zati The isomorphism $f$ is optimal and the Euler classes $e_{1}$
and $e_{2}$ are proportional.

\zati There exists a perverse isomorphism 
$\mbox{\boldmath $F$} \colon \lau{\mbox{\boldmath $\IH$}}{}{}{X_{2}}\to 
\lau{\mbox{\boldmath $\IH$}}{}{}{X_{1}}$
verifying 
$
\mbox{\boldmath $F$} \rondp \mbox{\boldmath $\pi_{2}$} = \mbox{\boldmath
$\pi_{1}$} \rondp \mbox{\boldmath $f$} $.  \ep 

\pro We proceed in two steps.

\medskip

\fbox{$(a) \Rightarrow (b)$} Since the isomorphism $f$ is optimal then
$\per{x}_1 = \per{x}_2$ and we will denote by $\per{x}$ this
perversity.  Since $f^{*}e_{2} = f^{*}[\epsilon_2] = \lambda \cdot e_{1} =
\lambda \cdot [\epsilon_1] $, with $\lambda \in \R\menos\{ 0 \}$, then there exists $\gamma \in
\lau{\Om}{1}{\per{e}}{B_{2}}$ with $f^{*}\epsilon_2 = \lambda \cdot
\epsilon_1 - d(f^{*}\gamma)$.
For each perversity $\per{p}$ we define
$\ib{F}{\per{p}}\colon \lau{I\Omega}{*}{\per{p}} {X_2}\TO 
\lau{I\Omega}{*}{\per{p}} {X_1} $
by
$$
\ib{F}{\per{p}}(\alpha,\beta)=(f^{*}(\alpha -\beta \wedge \gamma),
f^{*}(\lambda \cdot \beta)).
$$
The map $\ib{F}{\per{p}}$ is a well defined differential graded
morphism.  Let us see that. For each $(\alpha,\beta)
\in\lau{I\Omega}{*}{\per{p}}{X_2} $ and for each $S \in
\mathcal{S}_{X_{1}}^{^{sing}}$ we have

\begin{itemize}

\item[-] $f^{*}(\alpha -\beta \wedge \gamma )\in
\lau{\Pi}{*}{}{B_{1}}$.

\item[-] $f^{*}(\lambda \cdot \beta)
\in \lau{\Om}{*-1}{\per{p}-\per{x}}{B_{1}}$.

\item[-] $||f^{*}(\alpha -\beta \wedge \gamma )||_{\pi(S)} 
=
||\alpha -\beta \wedge \gamma ||_{\pi(f(S))}\leq 
\max \left( ||\alpha ||_{\pi(f(S))} , ||\beta ||_{\pi(f(S))}+ || \gamma ||_{\pi(f(S))}
\right) \\ \leq \max \left( \per{p}(S), \per{p}(S) -\per{x}(S) + ||
\gamma ||_{\pi(f(S))} \right) \leq\per{p}(S) $ since $|| \gamma ||_{\pi(f(S))} \leq
\per{x}(S) $.

\item[-] $||f^{*}d(\alpha -\beta \wedge \gamma) +(-1)^{|\beta|} f^{*}(\lambda \cdot \beta)
\wedge \epsilon_1 ||_{\pi(S)}
= 
||f^{*}d\alpha -f^{*}(d\beta \wedge \gamma)
-(-1)^{|\beta|}f^{*}(\beta \wedge d\gamma )+
(-1)^{|\beta|} f^{*} (\beta \wedge \epsilon_{2}) +
(-1)^{|\beta|}f^{*}(\beta \wedge d\gamma)||_{\pi(S)}
=
||f^{*}\left(d\alpha +(-1)^{|\beta|} \beta \wedge \epsilon_{2}\right)
-f^{*}(d\beta \wedge \gamma) ||_{\pi(S)} 
\leq \max ( ||d\alpha +(-1)^{|\beta|} \beta \wedge
\epsilon_{2}||_{\pi(f(S))}, ||d\beta \wedge \gamma)||_{\pi(f(S))} )
\leq
\per{p}(S).$

\item[-] $D_1 \ib{F}{\per{p}} (\alpha,\beta) = 
(f^{*}d(\alpha -\beta \wedge \gamma) +(-1)^{|\beta|} f^{*}(\lambda
\cdot \beta) \wedge \epsilon_1, f^{*}(\lambda \cdot  d\beta)) = 
(f^{*}\left(d\alpha +(-1)^{|\beta|} \beta \wedge \epsilon_{2}\right)
-f^{*}(d\beta \wedge \gamma),f^{*}(\lambda \cdot  d\beta))
= \ib{F}{\per{p}} (d\alpha + (-1)^{|\beta|}
\beta \wedge \epsilon_2,d\beta) = \ib{F}{\per{p}} D_2 (\alpha,\beta). 
$
\end{itemize}

\smallskip

\nt The family $\mbox{\boldmath $F$} = \{ \ib{F}{\per{p}} \} \colon
\lau{\mbox{\boldmath $I\Om$}}{}{}{X_{2}}\to \lau{\mbox{\boldmath
$I\Om$}}{}{}{X_{1}}$ is a perverse morphism since:

   \begin{itemize}
       \item[\refp{bat}] A straightforward calculation.  
       \item[\refp{bi}]
   Consider $ (\alpha,\beta) \in \lau{I\Om}{*}{\per{p}}{X_2}$ and
   $(\alpha',\beta')\in \lau{I\Om}{*}{\per{p}'}{X_2}$.  Then

$\ib{F}{\per{p}+ \per{p}'} ((\alpha,\beta) \wedge (\alpha',\beta')) =
\ib{F}{\per{p}+ \per{p}'} (( \alpha \wedge \alpha' ,
(-1)^{|\alpha'|}\beta\wedge \alpha' + \alpha \wedge \beta' )) = (
f^{*}(\alpha \wedge \alpha' - (-1)^{|\alpha'|}\beta\wedge \alpha' \wedge
\gamma - \alpha \wedge \beta' \wedge \gamma),
f^{*}((-1)^{|\alpha'|}\lambda \cdot \beta\wedge \alpha' + \lambda
\cdot \alpha \wedge \beta')) = ( f^{*}(\alpha -\beta \wedge \gamma) ,
f^{*}(\lambda \cdot \beta )) \wedge ( f^{*}(\alpha' -\beta' \wedge \gamma)
,f^{*}(\lambda\beta' ) )= \ib{F}{\per{p}}(\alpha,\beta) \wedge
\ib{F}{\per{p}'}(\alpha',\beta').  $
\end{itemize} 
In fact, the perverse morphism $\mbox{\boldmath $F$} $ is a perverse isomorphism,
the inverse is given by $\mbox{\boldmath $F$}^{-1} 
=
\{ \hiru{F}{-1}{\per{p}}\},
$
where 
$\hiru{F}{-1}{\per{p}}(\alpha,\beta)=(f^{-*}\alpha + \lambda^{-1}
\cdot f^{-*}\beta \wedge \gamma, \lambda^{-1}\cdot f^{-*}\beta).  $
We conclude that the induced operator 
$\mbox{\boldmath $F$} \colon
\lau{\mbox{\boldmath $\IH$}}{}{}{X_{2}}\to \lau{\mbox{\boldmath
$\IH$}}{}{}{X_{1}}
$
is a perverse isomorphism.  
Finally, the equality $ \mbox{\boldmath $F$} \rondp \mbox{\boldmath
$\pi_{2}$} = \mbox{\boldmath $\pi_{1}$} \rondp \mbox{\boldmath $f$} $
comes from
$$
 \ib{F}{\per{p}} ( \ib{(\pi_2)}{\per{p}}(\alpha))=
\ib{F}{\per{p}} (\alpha , 0) = (f^{*}\alpha , 0)
=\ib{(\pi_1)}{\per{p}}(f^{*}\alpha) = 
\ib{(\pi_1)}{\per{p}}(\ib{f}{\per{p}}(\alpha)), $$
where $\per{p}$ is a perversity and $\alpha \in
\lau{I\Om}{*}{\per{p}}{B_{2}}$.  

\medskip

\fbox{$(b) \Rightarrow (a)$}
Write
$\mbox{\boldmath $f$} = \left\{
\ib{f}{\per{p} }\colon \lau{\IH}{*}{\per{p}}{B_2} \TO
\lau{\IH}{*}{\per{p}}{B_1}\right\}$ and $\mbox{\boldmath $F$} =
\left\{ \ib{F}{\per{p} }\colon \lau{\IH}{*}{\per{p}}{X_2} \TO
\lau{\IH}{*}{\per{p}}{X_1}\right\}$.  Consider now the Gysin sequences
associated to the action $\Phi_1$ and $\Phi_2$.  The two Gysin terms
are written {$_{_1}\mathcal{G}$} and {$_{_2}\mathcal{G}$}
respectively.  Since $\ib{F}{\per{e}_{2}} \circ \ib{(\pi_2)}{\per{e}_{2}} =
\ib{(\pi_1)}{\per{e}_{2}} \circ \ib{f}{\per{e}_{2}}$ we can construct a
commutative diagram \be
\label{diag}
\begin{CD}
\lau{\IH}{1}{\per{e}_2}{B_{2}} 
@>\ib{(\pi_2)}{\per{e}_2}>>\lau{\IH}{1}{\per{e}_2}{X_2} 
@>\ib{(\oint_{2})}{\per{e}_2}>> 
\coho{H}{0}{\lau{_{_2}\mathcal{G}}{*}{\per{e}_2}{B_{2}}} 
@>\ib{(\eub_2)}{\per{e}_2}>> 
\lau{\IH}{2}{\per{e}_2}{B_{2}} 
@>\ib{(\pi_2)}{\per{e}_2}>>\lau{\IH}{2}{\per{e}_2}{X_2} \\
    @V \ib{f}{\per{e}_2} VV @V\ib{F}{\per{e}_2} VV @V\ell VV     
@V \ib{f}{\per{e}_2} VV @V\ib{F}{\per{e}_2} VV \\
\lau{\IH}{1}{\per{e}_2}{B_{1}} 
@>\ib{(\pi_1)}{\per{e}_2}>>\lau{\IH}{1}{\per{e}_2}{X_1} 
@>\ib{(\oint_{1})}{\per{e}_2}>>
\coho{H}{0}{\lau{_{_1}\mathcal{G}}{*}{\per{e}_2}{B_{1}}}
@>\ib{(\eub_1)}{\per{e}_2}>> \lau{\IH}{2}{\per{e}_2}{B_{1}}
@>\ib{(\pi_1)}{\per{e}_2}>>\lau{\IH}{2}{\per{e}_2}{X_1} ,\\
    \end{CD}
    \ee
where $\ell \colon \coho{H}{0}{\lau{_{_2}\mathcal{G}}{*}{\per{e}_2}{B_{2}}} 
\to 
\coho{H}{0}{\lau{_{_1}\mathcal{G}}{*}{\per{e}_2}{B_{1}}}$
 is an isomorphism. From  \refp{g} we get that
 $\coho{H}{0}{\lau{_{_2}\mathcal{G}}{*}{\per{e}_2}{B_{2}}} $ is $\R$
 (the constant functions) and therefore $\ell$ is the multiplication by a number
 $\lambda \in \R \backslash \{ 0\}$.  We prove (a) in two steps.  \smallskip

{\em 1. If the isomorphism $f$ is optimal then the Euler classes $e_{1}$
and $e_{2}$ are proportional}.
We have $\per{e}_1 = \per{e}_2 = \per{e}$. 
The formula \refp{g} and the diagram \refp{diag} give $$ \lambda \cdot e_1 = 
\lambda \cdot \ib{(\eub_1)}{\per{e}}(1)  = \ib{f}{\per{e}} \left(
\ib{(\eub_2)}{\per{e}}(1)\right) = \ib{f}{\per{e}} ( e_2).$$

\smallskip

{\em 2. The isomorphism $f$ is optimal}. It suffices to prove 
that $\per{e}_1 (\pi_{1}(S))= \per{e}_2(f(\pi_{1}(S)))$ for each $S \in {\cal
S}^{^{sing}}_{X_{1}}$.
Since
 $\lau{H}{0}{}{\lau{_{_1}\mathcal{G}}{*}{\per{e}_2}{B_{1}}}= \R$ then
 $1 \in \lau{_{_1}\mathcal{G}}{*}{\per{e}_2}{B_{1}}$  and we get that 
$\per{e}_2 -\per{x}_1 \geq 0$. So,
$\per{e}_1(\pi_{1}(S)) =0$ if $\per{e}_2(f(\pi_{1}(S)))= 0.$
By symmetry :
$
\per{e}_1 (\pi_{1}(S))= 0\Longleftrightarrow \per{e}_2(f(\pi_{1}(S)))=0.
$

The fixed strata are the same for both actions. If the perverse strata are 
different, then we can find a fixed stratum $S$ with $\per{e}_1(\pi_{1}(S) )\ne 
\per{e}_2(f(\pi_{1}(S)))$ and $\per{e}_1 (\pi_{1}(S'))= \per{e}_2(f(\pi_{1}(S')))$ 
for each singular stratum $S'$ with $S \preceq S'$.  In particular, the fixed
strata and the perverse strata are the same on $L_S$.  We have proved
that the Euler classes of the actions $\Phi_{1,L_S} \colon \sbat
\times L_S \to L_S$ and $\Phi_{2,L_S} \colon \sbat \times L_S \to L_S$
are proportional trough a non-vanishing factor.  So, they vanish or
not simultaneously.  This would give $\per{e}_1(\pi_{1}(S)) =
\per{e}_2(f(\pi_{1}(S)))$ (cf.  1.3).  Contradiction. \qed

\prgg {\bf Remark}.  The connectedness and the normality of $X_{1}$ and
$X_{2}$ have only been used in the proof of $(b) \Rightarrow (a)$. 

The
first result of this work shows how the Euler class of the action
determines the intersection cohomology algebra of the unfolded
pseudomanifold $X$.

\bc
\label{uno}
Consider two modelled actions $\Phi_1 \colon \sbat \times X_1 \to
X_1$ and $\Phi_2 \colon \sbat \times X _2\to X_2$ having a common
orbit space.  If the Euler classes $e_{1}$ and $e_{2}$ are
proportional then intersection cohomology algebra of $X_{1}$ and
$X_{2}$ are isomorphic.  \ec

The second result of this work shows how the Euler class of the action
determines the real homotopy type of the stratified unfolded $X$.

\bc 
\label{dos}
Let $X_{1}$, $X_{2}$ be two connected normal unfolded pseudomanifolds. Consider two modelled actions $\Phi_1 \colon \sbat \times X_1
\to X_1$ and $\Phi_2 \colon \sbat \times X _2\to X_2$ having a common
orbit space.  If the two Euler classes $e_{1}$ and $e_{2}$ are
proportional than the real homotopy type of $X_1$ and $X_2$ are the
same.  \ec

\pro The real homotopy type of $X_k$ is determined by the dgca 
$\lau{\Om}{*}{\per{0}}{X_k}$ for $k=1,2$ (cf. \cite{RS}).
The result comes from the following sequence of dgca quasi-isomorphisms:

$$
\lau{\Om}{*}{\per{0}}{X_2}
\stackrel{\ib{\Delta}{2,\per{0}}}{\OT \! \! \! - }
\lau{I\Om}{*}{\per{0}}{X_2}
\stackrel{\ib{F}{\per{0}}}{\TO}
\lau{I\Om}{*}{\per{0}}{X_1}
\stackrel{\ib{\Delta}{1,\per{0}}}{- \!  \! \!\TO}
\lau{\Om}{*}{\per{0}}{X_1}
$$
(cf. \refp{delta}, Proposition \ref{comp}).
\qed

Inspired by the notion of real homotopy type we can define the perverse
real homotopy type of an unfolded pseudomanifold in the following way.
Two unfolded pseudomanifolds $X_{1}$ and $X_{2}$ have the same {\em perverse
real homotopy type} if there exists a finite family of perverse
quasi-isomorphisms
$$
X_{1} \leftarrow \bullet \rightarrow \cdots \leftarrow\bullet \to X_{2}.
$$
Here, a perverse
quasi-isomorphism is a perverse isomorphism inducing an isomorphism in
cohomology.  Notice that, in the Proposition \ref{comp}, we have proved
in fact the following result:

\bc 
\label{tres}
Consider two modelled actions $\Phi_1 \colon \sbat \times X_1
\to X_1$ and $\Phi_2 \colon \sbat \times X _2\to X_2$ having a common
orbit space.  If the two Euler classes $e_{1}$ and $e_{2}$ are
proportional then the perverse real homotopy type of $X_1$ and $X_2$ are the
same.  \ec

\end{document}